\documentclass{article}
\usepackage{amsthm,amssymb,amsmath}
\usepackage{color}
\usepackage{epsfig,verbatim}

\newtheorem{thm}{Theorem}
\newtheorem{prop}[thm]{Proposition}
\newtheorem{defn}[thm]{Definition}
\newtheorem{lem}[thm]{Lemma}

\theoremstyle{remark}
\newtheorem*{rem}{Remark}

\title{The Calabi functional on a ruled surface}
\author{G\'abor Sz\'ekelyhidi}
\date{}

\begin{document}

\maketitle
\begin{abstract}
	We study the Calabi functional on a ruled surface over a genus two
	curve. For polarisations which do not admit an extremal metric we
	describe the behaviour of a minimising sequence splitting the
	manifold into pieces. We also show that the Calabi flow
	starting from a metric with suitable symmetry gives
	such a minimising sequence.  
\end{abstract}

\section{Introduction}

In~\cite{Cal82} Calabi introduced the problem of minimising the $L^2$-norm of
the scalar curvature (this is called the \emph{Calabi functional})
over metrics in a fixed K\"ahler class on a compact
K\"ahler manifold. A critical point of the Calabi functional is called an
extremal metric. The Euler-Lagrange equation is that the gradient of the scalar
curvature is a holomorphic vector field.
It is known that extremal metrics in fact minimise the Calabi
functional (see~\cite{Hwang95}, \cite{Chen06}, \cite{Don05}). 
Recently much progress has been made in understanding
when extremal metrics exist, at least on a conjectural level. K\"ahler-Einstein
metrics are a special case and when the first Chern class of the manifold is
positive (the manifold is called Fano in this case), 
Yau conjectured that the existence of
K\"ahler-Einstein metrics is 
related to the stability of the manifold in the sense of geometric invariant
theory. In the case of negative or zero first Chern class Yau~\cite{Yau78}
and Aubin~\cite{Aub78} have shown
that K\"ahler-Einstein metrics always exist, answering a conjecture of Calabi.
Tian~\cite{Tian97} made significant progress towards understanding the Fano
case, solving it completely in the case of surfaces in~\cite{Tian90}.
Donaldson~\cite{Don97} 
showed that the scalar curvature can be interpreted as a moment map (this was
also observed by Fujiki~\cite{Fuj92}) 
and this enabled extending the conjectures about the
existence of 
K\"ahler-Einstein metrics to more general constant scalar curvature and extremal
metrics (see~\cite{Don02}, \cite{Mab04_1}, \cite{GSz04}). 

In this paper we look at what we can say about minimising 
the Calabi functional in a K\"ahler class which admits no
extremal metric, concentrating on a concrete example. Let $\Sigma$ be a
genus 2
curve and $\mathcal{M}$ a degree -1 line bundle on it. We consider the ruled
surface $X=\mathbf{P}(\mathcal{M}\oplus\mathcal{O})$ with a family of
polarisations $L_m = C + mS_\infty$, where $C$ is the class of a fibre,
$S_\infty$ is the infinity section (with self-intersection 1), and $m>0$. 
Technically we should take $m$ to
be rational, especially when discussing test-configurations, but by an
approximation and continuity argument
we can take $m$ to be real. The aim is to study the problem
of minimising the Calabi functional in these K\"ahler classes. Our main result
is the following. 

\begin{thm}\label{thm:main} 
  There exist constants $k_1\simeq 18.9,k_2\simeq 5.03$, such that
  \begin{enumerate} 
    \item If $0<m<k_1$ then $X$ admits an extremal
      metric (this is due to T{\o}nnesen-Friedman~\cite{TF97}).  
    \item If
      $k_1\leqslant m\leqslant k_2(k_2+2)$ then there exists a minimising sequence of
      metrics which breaks $X$ into two pieces and converges to complete
      extremal metrics on both.  
    \item If $m>k_2(k_2+2)$ then there exists a minimising sequence of
      metrics which breaks $X$ into three pieces. It converges to complete
      extremal metrics on two of these and the third degenerates into a fibration of
      infinitely long and infinitely thin cylinders.  
  \end{enumerate}
\end{thm}

\begin{comment}
These decompositions are analogous to the Harder-Narasimhan filtration of an
unstable vector bundle. The extremal metrics on them correspond to the
Hermitian-Einstein metrics on the successive quotients of the Harder-Narasimhan
filtration. See~\cite{Don02} for the conjectural decomposition
of unstable toric varieties. 
\end{comment}

To construct metrics on our ruled surface, we use the momentum construction
given in Hwang-Singer~\cite{HS02}. This construction has been used repeatedly in
the past to
find special metrics on ruled manifolds, in particular extremal metrics.
See~\cite{ACGT3} for a unified treatment of these constructions or~\cite{HS02} for
a historical overview and more references. The momentum construction 
allows us to construct circle invariant
metrics from functions on an interval and it gives a convenient expression
for the scalar curvature. More precisely, let 
$\phi:[0,m]\to\mathbf{R}$ be a smooth function, positive on the
interior $(0,m)$, vanishing at the endpoints, and such that $\phi'(0)=2$,
$\phi'(m)=-2$. The momentum construction gives a metric $\omega_\phi$ in
the K\"ahler class $L_m$, with scalar curvature
\[ S(\omega_\phi) =  \frac{-2}{1+\tau}-\frac{1}{2(1+\tau)}
  \big[ (1+\tau)\phi\big]^{\prime\prime}. 
\]
Here $\tau$ is the moment map for the $S^1$-action on the fibres and working
with this coordinate is the central idea of the momentum construction.  
We will recall this construction in
Section~\ref{sec:momentum}. Of particular importance to us is the fact that we
can consider momentum profiles which vanish on a subset of $(0,m)$. These
correspond to degenerate metrics and they arise as the limits of the minimising
sequences in Theorem~\ref{thm:main}.

In Section~\ref{sec:minimising} 
we consider the problem of directly minimising the Calabi
functional on the set of metrics obtained by the momentum construction. Since
the $L^2$-norm of the scalar curvature is equivalent to
the $H^2$-norm of the momentum
profiles, this is straight forward. We find that the Euler-Lagrange
equation for a minimiser $\phi$ is $\phi S(\phi)''=0$ and $S(\phi)''$ must
be a negative distribution, ie. $S(\phi)$ is concave. We show that a
unique minimiser
exists in each K\"ahler class and its momentum profile is in $C^2$.
Note that $S(\phi)''=0$ is the equation for $\phi$ to define an extremal metric. 

In Section~\ref{sec:explicit} we
explicitly construct the minimisers, which can be degenerate in the sense that
the momentum profiles can vanish on a subset of $(0,m)$. Here we will see the three
different kinds of behaviour stated in Theorem~\ref{thm:main}. 
In Section~\ref{sec:test-configs} we construct test-configurations for $X$ and
calculate their Futaki invariants. This will clarify the role of the
concavity of $S(\phi)$ for minimisers of the Calabi functional. In fact,
rational, piecewise-linear convex functions on $[0,m]$ give
test-configurations essentially by the construction in~\cite{Don02} as
generalised to bundles of toric varieties in~\cite{GSzThesis}. We
can approximate $-S(\phi)$ by such functions, and 
Donaldson's theorem on lower bounds for
the Calabi functional in~\cite{Don05} shows that $\omega_\phi$ actually 
achieves the infimum of the Calabi functional on
the whole K\"ahler class, not just the metrics arising from the momentum
construction. This will complete the proof of Theorem~\ref{thm:main}.

An alternative approach to minimising the Calabi functional is using the Calabi
flow introduced in~\cite{Cal82}. This is the flow which deforms the K\"ahler
potential in the direction of the scalar curvature. It is
expected (see~\cite{Don02}, \cite{Don04_1}) that the Calabi flow should
minimise the Calabi functional and if there is no extremal metric in a
given K\"ahler class, then it should break up the manifold into pieces
which admit complete extremal metrics or collapse in some way.  
In Sections~\ref{sec:calabiflow} and \ref{sec:longtimeexist} we will
verify this, showing

\begin{thm} If the initial metric is given by the momentum construction then the
	Calabi flow exists for all time and converges
	to the minimiser of the Calabi functional. 
\end{thm}

The Calabi flow on ruled manifolds has been previously studied
in~\cite{Gua05}, where the long time existence and convergence is 
proved for the K\"ahler classes
which admit an extremal metric. We use similar techniques, the main
difference being that we introduce some variants of the Mabuchi
functional when no extremal metric exists.
In particular in the unstable case where $k_1 \leqslant m \leqslant
k_2(k_2+2)$ we define a functional which decreases along the Calabi
flow, is bounded below, and whose derivative is given by the difference
between the Calabi functional and its infimum. This leads to the
convergence result. The case $m > k_2(k_2+2)$ is more delicate since the
analogous Mabuchi-type functional is not bounded from below. Nevertheless
it has at worst logarithmic decay along the Calabi flow and this is
enough to show that the flow minimises the Calabi functional. 
This is discussed in Section~\ref{sec:calabiflow}.

Note that throughout the paper we have ignored factors of $2\pi$, for example in
the definition of the Calabi functional. Also we normalise the Futaki invariant
slightly differently from usual in Section~\ref{sec:test-configs}. Hopefully
this will lead to no confusion.

\subsubsection*{Acknowledgements}
Part of this work has appeared in the author's PhD thesis~\cite{GSzThesis}. I
would like to thank my supervisor Simon Donaldson for his encouragement and
for sharing his insights. 

\section{Metrics on the ruled surface} \label{sec:momentum}
In this section we describe the momentum construction for metrics on the
ruled surface (see Hwang-Singer~\cite{HS02}). 
Let $X$ be the ruled surface as above, so that
$X=\mathbf{P}(\mathcal{M}\oplus\mathcal{O})\to\Sigma$, where $\Sigma$ is a
genus 2 curve, and $\mathcal{M}$ is a degree -1 line bundle over
$\Sigma$. Let $\omega_\Sigma$ be a metric on $\Sigma$ with area $2\pi$
and constant scalar curvature $-2$. Also, let $h$ be a Hermitian metric on
$\mathcal{M}$ with curvature form $i\omega_\Sigma$. We consider metrics
on the total space of $\mathcal{M}$ of the form
\[  \omega = p^*\omega_{\Sigma} + 2i\partial\bar{\partial} f(s), \]
where $p:\mathcal{M}\to\Sigma$ is the projection map,
$s=\frac{1}{2}\log|z|^2_h$ is the logarithm of the fibrewise norm and
$f(s)$ is a suitable strictly convex function that makes $\omega$
positive definite.  
The point of the momentum construction is the change of coordinate from
$s$ to $\tau=f^\prime(s)$. The metric $\omega$ is invariant under the
$U(1)$-action on $\mathcal{M}$, and $\tau$ is just the moment map for
this action. Let $I\subset\mathbf{R}$ be the image of $\tau$, and let
$F:I\to\mathbf{R}$ be the Legendre transform of $f$. By definition this
means that
\[ f(s) + F(\tau) = s\tau, \]
and $F$ is a strictly convex function. The \emph{momentum profile} is
defined to be the function 
\[ \phi(\tau) = \frac{1}{F^{\prime\prime}(\tau)}.\]
We have the following relations:
\[ s = F^\prime(\tau),\quad \frac{ds}{d\tau}=F^{\prime\prime}(\tau),
\quad
\phi(\tau) = f^{\prime\prime}(s). \]

\subsubsection*{The metric in local coordinates}
Let us now see what the metric $\omega$ looks like in local coordinates.
Choose a local coordinate $z$ on $\Sigma$ and a fibre coordinate $w$ for
$\mathcal{M}$. The fibrewise norm is given by $|(z,w)|^2_h = |w|^2h(z)$
for some positive function $h$, so that 
\[ s = \frac{1}{2}\log |w|^2 + \frac{1}{2}\log h(z). \]
We can choose the local trivialisation of $\mathcal{M}$ in such a way that at 
a point $(z_0,w_0)$ we have
$d\log h(z)=0$. We can then compute at the point $(z_0,w_0)$
\[ 
\begin{aligned}
  2i\partial\bar{\partial} f(s) &= if^\prime(s)\partial\bar{\partial}
\log h(z) + f^{\prime\prime}(s)\frac{i\,dw\wedge d\bar{w}}{
2|w|^2} \\
&= \tau p^*\omega_\Sigma + \phi(\tau)\frac{i\,dw\wedge d\bar{w}}{
2|w|^2}.
\end{aligned}
\]
The metric at the point $(z_0,w_0)$ is therefore given by
\begin{equation}\label{eq:momentmetric}
  \omega = (1+\tau)p^*\omega_\Sigma + \phi(\tau)\frac{i\,dw\wedge
d\bar{w}}{2|w|^2}.
\end{equation}
In order to compute the scalar curvature of $\omega$, note that the
determinant of the metric $g$ correponding to $\omega$ is
\[ \det(g) = \frac{1}{|w|^2}(1+\tau)\phi(\tau)\det(g_\Sigma), \]
which is valid for all points, not just $(z_0,w_0)$. The Ricci form at
$(z_0,w_0)$ is given by
\[ 
\begin{aligned}
  \rho &= -i\partial\bar{\partial}\log\det g \\
  &= p^*\rho_\Sigma - \frac{\big[
  (1+\tau)\phi\big]^{\prime}}{2(1+\tau)}
  p^*\omega_\Sigma-\frac{\phi}{2}\cdot \frac{
  (1+\tau)\big[ (1+\tau)\phi\big]^{\prime\prime} -
  \big[ (1+\tau)\phi\big]^\prime}{(1+\tau)^2}\cdot\frac{i\, dw\wedge
  d\overline{w}}{\vert w\vert^2},
\end{aligned}
\]
where the derivatives are all with respect to $\tau$ (note
that $d/ds = \phi(\tau) d/d\tau$) and $\rho_\Sigma$ is the Ricci form of
the metric $\omega_\Sigma$. 
Taking the trace of this, we find that the scalar curvature $S(\omega)$
is given by
\begin{equation}\label{eq:scalarcurv}
  S(\omega) = \frac{-2}{1+\tau}-\frac{1}{2(1+\tau)}
  \big[ (1+\tau)\phi\big]^{\prime\prime}.
\end{equation}

In~\cite{HS02} the extendability of the metrics to the projective completion of
$\mathcal{M}$ is studied. The proposition we need is the following. 

\begin{prop}[see~\cite{HS02}] \label{prop:momentumprofile}
  For some $m>0$ let
  $\phi:[0,m]\to\mathbf{R}$ be a smooth function such that $\phi$ is
  positive on $(0,m)$, and
  \begin{equation}\label{eq:bdrycond}
	  \phi(0)=\phi(m)=0,\qquad \phi^\prime(0)=2,\quad \phi^\prime(m)=-2.
  \end{equation}
  Then the momentum construction defines a smooth metric $\omega_\phi$
  on $X$ in the K\"ahler class $C+mS_\infty$, with scalar curvature
  $S(\omega)(\tau)$ given by Equation~\ref{eq:scalarcurv}. Here $C$ is
  the class of a fibre, and $S_\infty$ the infinity section.

  If instead $\phi$ satisfies the boundary conditions
   \[ \phi(0)=\phi(m)=0,\qquad \phi^\prime(0)=0,\quad \phi^\prime(m)=-2,
  \]
  and $\phi(\tau)\leqslant O(\tau^2)$ for small $\tau$, 
  then the momentum construction gives a complete metric with finite
  volume on the complement of the zero section in $X$. Similarly if
  $\phi^\prime(0)=2$ and $\phi^\prime(m)=0$ then we obtain a complete
  metric on the complement of the infinity section. 

  The metrics are extremal, ie. their scalar curvature has holomorphic gradient, 
  when $S(\phi)''=0$. 
\end{prop}

Let us also note the definition
\begin{defn}\label{defn:momprof}
	A \emph{momentum profile} is a $C^2$ function $\phi : [0,m]\to
	\mathbf{R}$ which is positive on $(0,m)$ and satisfies the boundary
	conditions (\ref{eq:bdrycond}). A \emph{singular momentum profile} is
	the same except we only require it to be non-negaive instead of
	positive, ie. it can vanish on a subset of $(0,m)$. 
\end{defn}

Let us write $\Phi$ for the unique solution of $S(\Phi)''=0$ satisfying
the same boundary conditions as a momentum profile.  Then $\Phi$ is positive on
$(0,m)$ precisely when the polarisation admits an extremal metric. 
We define the Calabi
functional to be

\[\begin{aligned}
	Cal(\phi) &= \int_0^m (S(\phi)-S(\Phi))^2(1+\tau)\,d\tau \\ 
	&= \int_0^m \frac{1}{4(1+\tau)}\left[\left( (1+\tau)(\Phi-\phi)
	\right)''\right]^2\, d\tau.
\end{aligned} \] 

\noindent This differs from the $L^2$-norm of $S(\phi)$ by a constant, since
\[
\int_0^m (S(\phi)-S(\Phi)) S(\Phi)\, (1+\tau)d\tau = \int_0^m \frac{1}{2}
[(1+\tau)(\Phi-\phi)]'' S(\Phi)\, d\tau = 0, \]
integrating by parts, so
\[ Cal(\phi) = \int_0^m S(\phi)^2\, (1+\tau)d\tau - \int_0^m S(\Phi)^2\,
(1+\tau)d\tau. \]

Throughout the paper when we integrate a function over $X$ which only depends on
$\tau$ we will often use the volume form $(1+\tau)d\tau$. From the formula
(\ref{eq:momentmetric}) we see that this is a constant multiple of 
the integral with respect to the volume form $\omega^2$. Because
of the boundary conditions on $\phi$ the Poincar\'e inequality shows that the 
Calabi functional is equivalent to the
$H^2$-norm of $\phi$. 
This makes it easy to minimise the Calabi functional directly as we
do in the next section.

\section{Minimising the Calabi functional}\label{sec:minimising}
It is fairly simple to directly minimise the Calabi functional on the set of
metrics which are given by momentum profiles. 
We introduce the set of functions
\[ A = \left\{ \phi:[0,m]\to\mathbf{R}\,\left|\begin{aligned} &\phi\in H^2, 
  \ \phi\geqslant0 \text{ and }
  \phi\text{ satisfies the }\\ &\text{boundary conditions
	in Proposition~\ref{prop:momentumprofile}}\end{aligned}\right. \right\}, \]
and we want to minimise the Calabi functional on this space. Let us choose a
minimising sequence $\phi_k\in A$. We have a bound $\Vert\phi_k\Vert_{H^2}
\leqslant C\cdot Cal(\phi_k)$, so we can choose a subsequence converging weakly to
some $\phi\in H^2$. Weak convergence in $H^2$ implies convergence in $C^1$ so
the boundary conditions and non-negativity hold in the limit, ie. 
$\phi\in A$. Moreover $Cal$ is lower-semicontinuous because the
$H^2$-norm is, so $\phi$ is the required minimiser. 

\begin{prop}\label{prop:minimumprofile}
	The minimiser $\phi$ in $A$ satisfies 
	$\phi S(\phi)''=0$ and $S(\phi)''$ is
	a negative distribution. In particular $S(\phi)$ is continuous, so 
	$\phi\in C^2$. Conversely if $\psi S(\psi)''=0$ and $S(\psi)$ is
	concave, then $\psi=\phi$. 
\end{prop}
\begin{proof}
	The variation of $Cal$ at $\phi$ is given by
	\[ DCal_\phi(\widetilde{\phi}) = -\int_0^m (S(\phi)-S(\Phi))\left[
	(1+\tau)\widetilde{\phi}\right]''\,d\tau.\]
	We are considering variations inside $A$, so $\widetilde{\phi}$ and its
	first derivative vanishes at the endpoints. We can therefore integrate
	by parts, and find that 
	\[ -\int_0^m S(\phi)''\widetilde{\phi}(1+\tau)\, d\tau \geqslant 0\]
	for all $\widetilde{\phi}$ such that $\phi+\epsilon\widetilde{\phi}\in
	A$ for small enough $\epsilon$. We can choose $\widetilde{\phi}$ to be
	an arbitrary non-negative smooth function which vanishes along with its
	first derivative at the endpoints. This shows that $S(\phi)''$ is a
	negative distribution. On the open set where $\phi$ is positive
	we can choose $\widetilde{\phi}$ to be negative or positive, so it
	follows that $S(\phi)''=0$ at these points. Therefore $\phi S(\phi)''=0$
	on $(0,m)$. The continuity of $S(\phi)$ follows from it being concave,
	and this implies that $\phi\in C^2$. 

	The converse follows from the following computation.
        \[ \begin{aligned}
	    Cal(\psi) & \leqslant Cal(\psi) + 
	    \int_0^m (S(\phi)-S(\psi))^2(1+\tau)\, d\tau \\
	    &= Cal(\phi) + 2\int_0^m
	    (S(\psi)-S(\phi))S(\psi)(1+\tau)\, d\tau \\
	    &= Cal(\phi) + \int_0^m \big[ (1+\tau)\phi - (1+\tau)\psi
	    \big]'' S(\psi)\, d\tau \\
	    &= Cal(\phi) + \int_0^m \phi S(\psi)'' (1+\tau)\, d\tau \\
	    &\leqslant Cal(\phi).
	\end{aligned} \]
  Since $Cal(\phi)$ is minimal we must have equality, ie. 
  \[\int_0^m (S(\phi)-S(\psi))^2 (1+\tau)\, d\tau = 0.\]
  This implies that $S(\phi)=S(\psi)$, from which it follows that 
  $\phi=\psi$.
\end{proof}

\section{Explicit minimisers}\label{sec:explicit}
In this section we compute explicitly the minimisers of the Calabi functional
for all polarisations. For each $m$ we are looking for a singular momentum
profile (Definition~\ref{defn:momprof}) 
%a $C^2$ function
%$\phi:[0,m]\to\mathbf{R}$ which satisfies the boundary conditions of a momentum
%profile, is non-negative, and 
such that $S(\phi)''=0$ wherever $\phi$ does not vanish, and in addition 
$S(\phi)$ is concave. 

There are three cases to consider depending on the polarisation. 

\subsubsection*{Case 1. There exists an extremal metric, $m<k_1\simeq
18.889$}

\begin{figure}[htbp]
    \begin{center}
      \input{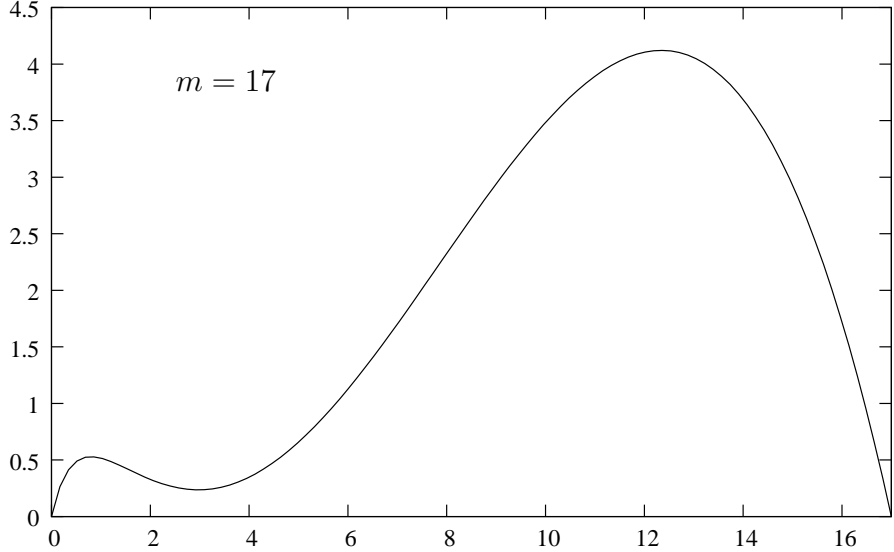}
      \caption{Momentum profile of an extremal metric on $X$ when
      $m=17$.}
    \label{fig:smooth}
    \end{center}
  \end{figure}

In this case we want to solve the equation $S(\phi)''=0$. By the
Formula (\ref{eq:scalarcurv}) for the scalar curvature, this is the ODE
\[ \frac{1}{2(1+\tau)}(-4-[(1+\tau)\phi]^{\prime\prime}) = A\tau + B,\]
for some constants $A, B$. Rearranging this and integrating twice we obtain
\begin{equation}\label{eq:ODE}
	(1+\tau)\phi = -\frac{A\tau^4}{6}-\frac{(A+B)\tau^3}{3}-B\tau^2-2\tau^2 +
	C\tau + D,
\end{equation}
where $C$ and $D$ are also constants. The boundary conditions on $\phi$ on the
interval $[0,m]$ give a
system of linear equations on $A,B,C,D$ which we can solve to obtain
\[\begin{aligned}
	\phi(\tau) = \frac{2\tau(m-\tau)}{m(m^2+6m+6)(1+\tau)}\big[& \tau^2(2m+2) +
	\tau(-m^2+4m+6)\\ & + m^2 + 6m+6\big].
\end{aligned}
\]
This will give a metric when it is positive on the interval $(0,m)$ which
happens if and only if the quadratic expression in square brackets is positive
on this interval. This is the case for $m<k_1$ where $k_1$ is the only positive
real roof of the quartic $m^4-16m^3-52m^2-48m-12$. Approximately
$k_1\simeq
18.889$, which is the result obtained by
T{\o}nessen-Friedman~\cite{TF97}. See Figure~\ref{fig:smooth} for a
graph of $\phi(\tau)$ for $m=17$. 

\subsubsection*{Case 2. $X$ breaks up into two pieces, $k_1\leqslant m\leqslant
k_2(k_2+2)\simeq 35.33$}

When $m\geqslant k_1$ we can no longer find a positive solution of $S(\phi)''=0$ on
the whole interval $[0,m]$ so we split the interval into two pieces 
$[0,c]$ and $[c,m]$. We would like to find $\phi$ which vanishes at $c$, but on the
intervals $(0,c)$ and $(c,m)$ we have $S(\phi)''=0$, and
$S(\phi)$ is concave on $[0,m]$. We first let $\phi_1$ be the solution
of the equation
\[\begin{gathered}
  S(\phi_1)''=0 \text{ on the interval } (0,c)\\
  \phi_1(0)=\phi_1(c)=0,\quad \phi_1'(0)=2,\quad \phi_1'(c)=0.
\end{gathered}\]
We obtain
\[
  \phi_1(\tau) = \frac{2\tau(c-\tau)^2}{c^2(c^2+6c+6)(1+\tau)}\big[
  \tau(-c^2+2c+3)+c^2+6c+6\big].
\]
This is positive on $(0,c)$ if the linear expression in square brackets is
positive on this interval. This happens for $c\leqslant k_2$ where $k_2$ is
the only positive real root of the cubic $c^2-3c^2-9c-6$. Approximately
$k_2\simeq 5.0275$.  The scalar curvature is given by 
\[ S(\phi_1) = \frac{12(c^2-2c-3)}{c^2(c^2+6c+6)}\tau -
\frac{6(2c^2-c-4)}{c(c^2+6c+6)}. \]

To deal with the interval $[c,m]$ we first solve the equation
\[\begin{gathered}
  S(\psi)''=0 \text{ on the interval } (0,d)\\
  \psi(0)=\psi(d)=0,\quad \psi'(0)=0,\quad \psi'(d)=-2.
\end{gathered}\]
for some constant $d$, and then shift the solution to $[c,m]$. The
solution on $[0,d]$ is given by
\[ 
  \psi(\tau) =
  \frac{2\tau^2(d-\tau)}{d^2(d^2+6d+6)(1+\tau)}\big[\tau(2d^2+4d+3)
  -d^3+3d^2+9d+6\big].
\]
As before, this is positive on $(0,d)$ if the linear term in square
brackets is
positive on this interval. This is the case for $d\leqslant k_2$, for the
same $k_2$ as above. The scalar curvature is given by
\[ S(\psi) = \frac{12(2d^2+4d+3)}{d^2(d^2+6d+6)}\tau - \frac{6(3d^2+5d
+2)}{d(d^2+6d+6)}. \]

Now note that if we define $\phi_2$ by
\[ \phi_2(\tau) = (c+1)\psi\left(\frac{\tau-c}{c+1}\right), \]
then $\phi_2$ solves the equation 
\[ \begin{gathered}
  S(\phi_2)''=0 \text{ on the interval } (c,(c+1)d+c) \\
  \phi_2(c)=\phi_2( (c+1)d+c )=0, \quad \phi_2'(c)=0,\quad \phi_2'(
  (c+1)d+c )=-2. 
\end{gathered}\]
The scalar curvature is given by 
\[ S(\phi_2)(\tau) = \frac{1}{c+1}\, S(\psi)\left(\frac{\tau-c}{c+1}
\right). \]

\begin{figure}[htbp]
    \begin{center}
      \input{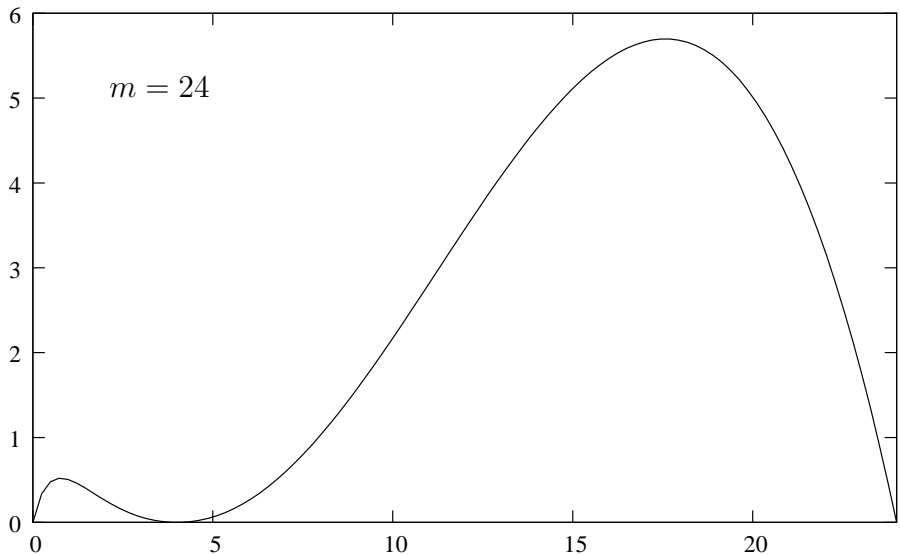}
      \caption{Momentum profile of the minimiser on $X$ when
      $m=24$. The manifold breaks into two pieces both of which are
      equipped with a complete extremal metric.}
    \label{fig:case2}
    \end{center}
  \end{figure}

We now define $\phi$ by
\[ \phi(\tau) = \begin{cases} \phi_1(\tau) \quad \tau\in[0,c],\\
  \phi_2(\tau) \quad \tau\in [c, (c+1)d + c].
\end{cases} \]
We can check that $S(\phi)$ will be continuous at $\tau=c$ precisely
when $c=d$. We also want $(c+1)d + c = m$, which implies that
$c=\sqrt{m+1}-1$. With these choices a simple computation shows that
$S(\phi)$ is concave for $m\geqslant k_1$ (note that it is linear for $m=k_1$,
and convex for $m<k_1$). Finally recall that the condition that $\phi$ is
non-negative means that $c\leqslant k_2$, which in turn implies $m
\leqslant k_2(k_2+2)$. See Figure~\ref{fig:case2} for a graph of $\phi$
for $m=24$. 

\subsubsection*{Case 3. $X$ breaks up into three pieces, $m >
k_2(k_2+2)$}

The previous construction no longer works for $m > k_2(k_2+2)$ so we
need to split the interval $[0,m]$ into three pieces. From the previous
case we have a solution $\phi_1$ to the equation
\[ \begin{gathered}
   S(\phi_1)''=0 \text{ on the interval } (0,k_1)\\
  \phi_1(0)=\phi_1(k_1)=0,\quad \phi_1'(0)=2,\quad \phi_1'(k_1)=0,
\end{gathered}\]
and also a solution $\phi_2$ to
\[ \begin{gathered}
   S(\phi_2)''=0 \text{ on the interval } (c,m)\\
  \phi_2(c)=\phi_2(m)=0,\quad \phi_2'(c)=0,\quad \phi_2'(m)=-2,
\end{gathered}\]
where the constant $c$ is defined by 
\begin{equation}\label{eq:defc}
  c = \frac{m+1}{k_2+1}-1.
\end{equation}
  
\begin{figure}[htbp]
    \begin{center}
      \input{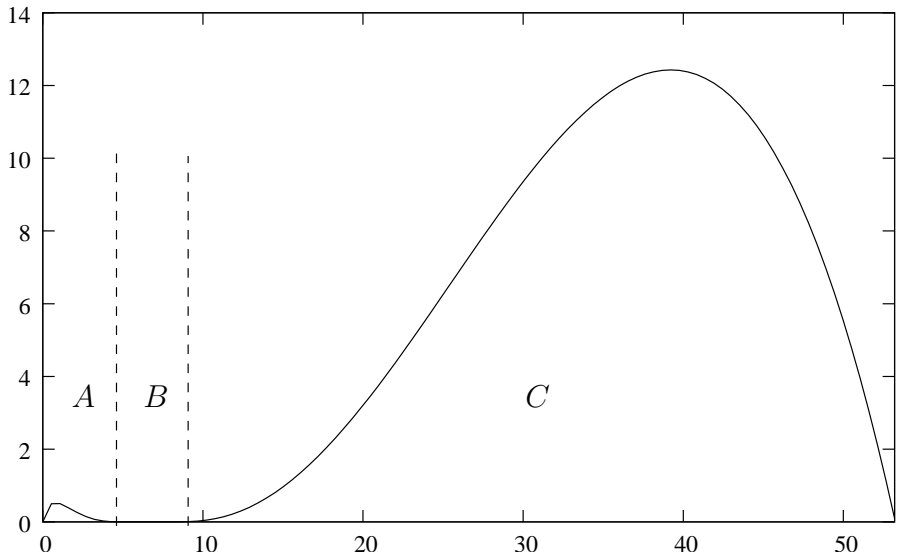}
      \caption{Momentum profile of the minimiser on $X$ when
      $m\simeq 53.2$. The manifold breaks into three pieces, two of
      which, $A$ and $C$,
      admit complete extremal metrics, and in the third, $B$, the
      $S^1$-orbits collapse.}
    \label{fig:case3}
    \end{center}
  \end{figure}

We define
\[ \phi(\tau) = \begin{cases} \phi_1(\tau)\quad &\tau\in[0,k_2] \\
  			0\quad &\tau\in[k_2,c]\\
			\phi_2(\tau)\quad &\tau\in[c,m].
\end{cases}\]
We can check that $c>k_2$ precisely when $m>k_2(k_2+2)$, and this choice
of $\phi$ satisfies that $\phi S(\phi)''=0$ and $S(\phi)$ is concave.
See Figure~\ref{fig:case3} for a graph of $\phi$ for $m\simeq 41.2$.

\subsubsection*{Conclusion}
For any $m$ one of the previous 3 cases will hold, so we can construct a
$\phi$ which satisfies the equation $\phi S(\phi)''=0$ and $S(\phi)$ is
concave. According to Proposition~\ref{prop:minimumprofile} this $\phi$ will
give the minimum of the Calabi functional on the space of singular
momentum profiles. In the next section we will show 
that they give the infimum of the Calabi
functional over all metrics in their K\"ahler class. This will complete
the proof of Theorem~\ref{thm:main}.

\section{Test-configurations}\label{sec:test-configs}
In the previous section we have found a (possibly degenerate) metric in each
K\"ahler class, which  minimises the Calabi functional on the set of 
metrics which come from the momentum construction. In this section we want to
show that these metrics minimise the Calabi functional on their entire K\"ahler
class. For this we use the theorem of Donaldson~\cite{Don05} which gives a lower
bound on the Calabi functional, given a destabilising test-configuration. We
will not give a detailed explanation of the test-configurations that we use, and
the computation of their Futaki invariants. For
more details see~\cite{GSzThesis} and \cite{Don02}.

\begin{prop}[Donaldson~\cite{Don05}]\label{prop:donaldson} Suppose there
  exists a test-configuration $\chi$ for a polarised variety $(X,L)$
  such that the Futaki invariant $F(\chi)$ is negative.  Then for any
  metric $\omega$ in the class $c_1(L)$ we have the inequality \[ \Vert
  S(\omega)-\hat{S}\Vert_{L^2} \geqslant \frac{-F(\chi_i)}{\Vert\chi_i\Vert}.\]
\end{prop}

The idea is to produce a sequence of test-configurations $\chi_i$ for which 
\[ \lim_{i\to\infty} \frac{-F(\chi_i)}{\Vert\chi_i\Vert} = \Vert S(\omega) -
\hat{S}\Vert_{L^2}, \]
where $\omega$ is the degenerate metric corresponding to the singular momentum
profile in each K\"ahler class that we have found in the previous
section. This will imply that this is the
infimum of the Calabi functional and $\omega$ minimises the Calabi
functional on its K\"ahler class. 

To obtain test-configurations we use the construction
in~\cite{GSzThesis} Section
4.1 (Theorem 4.1.2), which is an extension of the construction of
test-configurations for toric varieties by Donaldson~\cite{Don02} to bundles of
toric varieties. For the case of our ruled surface we obtain
\begin{prop} Given a rational, piecewise-linear, convex function
	$h:[0,m]\to\mathbf{R}$, there exists a test-configuration for
	$(X,L_m)$ with Futaki invariant given by 
	\begin{equation} \label{eq:convexfutaki}
		F(h) = h(0) + (1+m)h(m) -2\int_0^m h(\tau)\,d\tau -
		\hat{S}\int_0^m h(\tau) (1+\tau)\, d\tau, 
	\end{equation}
	and norm 
	\[ \Vert h\Vert^2 = \int_0^m (h(\tau)-\hat{h})^2(1+\tau)\,d\tau, \]
	where $\hat{h}$ is the average of $h$ with respect to the
	measure $(1+\tau)\,d\tau$.
\end{prop}

\begin{comment}
We briefly describe these test-configurations now. They are essentially
test-configurations for the $\mathbf{P}^1$ fibres. Given 
$h$ choose a large number integer $R$ such that $R-h(\tau)$ is positive on
$[0,m]$. Consider the polygon $Q\in\mathbf{R}^2$ defined by
\[ Q = \{(\tau, y)\,|\, y\leqslant h(\tau)\}. \]
By rescaling the lattice if necessary (this corresponds to taking a power of
$L_m$ as our polarisation) we can assume that $Q$ is an integral polytope, and
defines a polarised toric variety $(V_Q,L_Q)$. 
\end{comment}

To work with test-configurations we should restrict to polarisations $L_m$ 
with $m$ rational but an approximation argument
gives us the conclusion of Proposition~\ref{prop:donaldson} 
for any real $m$ as well.
Given a continuous convex function $h$ on $[0,m]$ which is not
necessarily rational or piecewise-linear, we still define the ``Futaki
invariant'' $F(h)$ of $h$ by Equation~\ref{eq:convexfutaki}. 

\begin{comment}
We very briefly indicate what kind of test-configurations these are. Given a
rational, piecewise-linear convex function $h$ on $[0,m]$, Donaldson's
construction gives a test-configuration for $\mathbf{P}^1$. Let us call the
total space of this test-configuration $\Chi$. 
\end{comment}

\begin{comment}
We have seen that circle invariant metrics on $X$ in the K\"ahler class
$L_m=C+mS_\infty$ correspond to smooth functions $\phi$ on the interval
$[0,m]$, positive on the interior, vanishing at $0,m$ and with
$\phi^\prime(0)=2, \phi^\prime(m)=-2$. Let us define a \emph{singular
momentum profile} to be a $C^2$ function on $[0,m]$ with the same
boundary conditions, but possibly vanishing on a subset of $(0,m)$. For
a singular momentum profile $\phi$ we still define the ``scalar
curvature'' of $\phi$ by Equation~\ref{eq:scalarcurv}.
\end{comment}

\begin{lem}\label{lem:futakiintegral}
  Let $\phi$ be a singular momentum profile, and
  $h:[0,m]\to\mathbf{R}$ a piecewise-smooth convex function. 
  Suppose that $h$ is linear on any interval on which $\phi$ does not
  vanish identically. Then 
  \[ F(h) = \int_0^m h(\tau)(S(\phi)-\hat{S})\,(1+\tau)d\tau.\]
\end{lem}

This result is analogous to the fact that the Futaki invariant of a
holomorphic vector field can be computed algebro-geometrically or
differential geometrically (see~\cite{Don02}).  Here if
$h$ is rational and piecewise-linear then it does not define a
holomorphic vector field but the result says that we
can still compute the Futaki invariant of the test-configuration it
induces with a differential geometric formula as long as we use a metric
which degenerates in a suitable way at points where $h$ is not linear.

\begin{proof} The proof is a simple integration by parts, using the
  formulas for $F(h)$ and $S(\phi)$.
\end{proof}

We can now complete the proof of Theorem~\ref{thm:main}.

\begin{proof}[Proof of Theorem~\ref{thm:main}]
  	What remains to be shown is that for each polarisation, the
	minimiser $\phi$ that we have constructed in the previous
	section minimises the Calabi functional over the whole K\"ahler
	class, not just over the set of metrics obtained from the
	momentum construction. Let $\phi$ be one of these minimisers.
	Since $-S(\phi)$ is convex, we can approximate it in the $C^0$-norm by
	a sequence of rational, piecewise-linear convex functions $h_i$. These
	define a sequence of test-configurations $\chi_i$ such that 
	\[ \lim_{i\to\infty}\frac{-F(\chi_i)}{\Vert\chi_i\Vert} = \frac{
	-F(-S(\phi)) }{ \Vert S(\phi) - \hat{S}\Vert_{L^2}}. \]
	If we let $h=-S(\phi)$, then $\phi$ and $h$ satisfy the conditions of
	Lemma~\ref{lem:futakiintegral} so that 
	\[ F(-S(\phi)) = -\int_0^m S(\phi)(S(\phi)(\tau)-\hat{S})(1+\tau)\, d\tau 
	= -\Vert S(\phi)-\hat{S}\Vert_{L^2}^2.\]
	Therefore
	\[\lim_{i\to\infty}\frac{-F(\chi_i)}{\Vert\chi_i\Vert} = \Vert
        S(\phi)-\hat{S}\Vert_{L^2},\]
	so that Proposition~\ref{prop:donaldson} now implies that this limit is
	the infimum of the Calabi functional on the K\"ahler class.
\end{proof}

\begin{comment}
\begin{rem}
Note that the test-configurations we defined are the worst destabilising
ones in the sense that they minimise the normalised Futaki invariant
$F(\chi)/\Vert\chi\Vert$. In this sense they are analogous to the
Harder-Narasimhan filtration of an unstable vector bundle
(see~\cite{BT05} for example). Geometrically the test-configurations we
obtain in the three cases are as follows. In the first case the
test-configuration is a product configuration induced by a holomorphic
vector field (the extremal vector field). In the second case it is
deformation to the normal cone of the zero (or infinity) section
(see~\cite{RT06} for the definition), but 
the $\mathbf{C}^*$-action on the total space of the test-configuration
is multiplied by a multiple of the product $\mathbf{C}^*$-action induced
by the extremal vector field. The central fibre in this case is two
copies of $X$ with a normal crossing along the zero section of one and
the infinity section of the other. In the third case approximating
$-S(\phi)$ by piecewise-linear convex functions
\end{rem}
\end{comment}

\section{The Calabi flow} \label{sec:calabiflow}
We have seen that in the case of a ruled surface it is fairly simple to minimise
the Calabi functional directly
over the set of metrics given by momentum profiles. It is
also interesting to see whether the Calabi flow converges to these minimisers.
In this section we will prove that this is the case. In~\cite{Gua05} Guan has
shown that on a ruled manifold
when an extremal metric exists, then starting from a metric given by the
momentum construction the Calabi flow exists for all
time and converges to the extremal metric exponentially fast. Our techniques are
similar to his, but we need to introduce some new functionals which are
modifications
of the Mabuchi functional more suited for studying the unstable
polarisations. 

We consider a family of metrics $\omega_s$ given by the momentum construction (see
Section~\ref{sec:momentum}), ie.
\[ \omega_t = p^*\omega_{\Sigma} + 2i\partial\overline{\partial} f_t(s), \]
for some family of suitably convex functions $f_t$. This path of metrics
satisfies the Calabi flow if
\[ \frac{\partial f_t}{\partial t} = S(\omega_t).\]
If we denote by $F_t$ the Legendre transforms of the $f_t$, then from the
definition of the Legendre transformation we find
\[ \frac{\partial F_t}{\partial t} = -\frac{\partial f_t}{\partial t}, \]
so that the path of momentum profiles $\phi_t = 1/F_t''$ satisfies
\[ \frac{\partial\phi_t}{\partial t} = \phi_t^2 S(\phi_t)'', \]
where $S(\phi_t)$ is given by Equation~\ref{eq:scalarcurv}.

It is known that the flow exists for a short time with any smooth initial metric
(see Chen-He~\cite{CH06}). 
Also, the Calabi functional is decreased under the flow:
\begin{lem}\label{lem:Calabidec}
	If $\phi$ is a solution to the Calabi flow, then
  \[
% \begin{aligned}
%    \frac{d\mathcal{M}(\phi)}{dt} & =-2\int_0^m
%    (S(\phi)-S(\Phi))^2(1+\tau)\, d\tau\leqslant 0, \\
    \frac{d\, Cal(\phi)}{dt} = -\int_0^m \phi^2\left( S(\phi)'' \right)^2
    (1+\tau)\,d\tau \leqslant 0.
%  \end{aligned}
  \]
  In particular the $H^2$ norm of $\phi_t$ is uniformly bounded along the flow. 
\end{lem}
\begin{proof}
  The result follows from the following computation of the variation.
  
  \begin{comment}
  \[\begin{aligned}
    \frac{d\mathcal{M}(\phi)}{dt} &= \int_0^m \left(-\Phi S(\phi)'' + \phi
    S(\phi)''\right)(1+\tau)\,d\tau \\
    &= \int_0^m (\phi-\Phi)(S(\phi)-S(\Phi))''(1+\tau)\,d\tau \\
    &= \int_0^m
    \left[(1+\tau)(\phi-\Phi)\right]''(S(\phi)-S(\Phi))\,d\tau\\
    &= -2\int_0^m (S(\phi)-S(\Phi))^2(1+\tau)\,d\tau.
  \end{aligned}
  \]
  For the second line note that $S(\Phi)''=0$, and in the next line we
  can integrate by parts twice because $\phi-\Phi$ and its first derivative
  vanish at the endpoints. 
\end{comment}
  
  \[ \begin{aligned}
    \frac{d\, Cal(\phi)}{dt} &= 2\int_0^m
    (S(\phi)-S(\Phi))\left(-\frac{1}{2(1+\tau)}\left[ (1+\tau)\phi^2
    S(\phi)'' \right]''\right)(1+\tau)\, d\tau \\
    &= - \int_0^m \phi^2\left(S(\phi)''\right)^2 (1+\tau)\,d\tau.
  \end{aligned}
  \]
  We can perform the integration by parts because $\phi^2$ and
  $(\phi^2)'$ vanish at the endpoints. Also recall that $S(\Phi)''=0$. 
\end{proof}

In Section~\ref{sec:longtimeexist} we will show that there is a solution to the
Calabi flow for all time for any polarisation. In this section we concentrate on
proving the following.

\begin{prop} If the flow exists for all time then the momentum profiles
  	converge in $H^2$ to the minimiser that we found in
	Section~\ref{sec:explicit}. 
\end{prop}
\begin{proof}
	Let us write $\Psi$ for the minimiser, so when $m < k_1$ then $\Psi$ is
	the momentum profile of an extremal metric, when $m\leqslant
	m\leqslant k_2(k_2+2)$
	then $\Psi$ vanishes at an interior point of $(0,m)$ and when
	$m > k_2(k_2+2)$ then $\Psi$ vanishes on an interval inside
	$(0,m)$. 

	Introduce the functional  
	\begin{equation}\label{eq:modmabuchi}
	  \mathcal{M}(\phi) = \int_0^m\left( \frac{\Psi}{\phi} +
	  \log\phi\right)(1+\tau)\,d\tau,
        \end{equation}
	defined on momentum profiles $\phi$. 
	When $m < k_1$ then in fact $\mathcal{M}$ is the modified Mabuchi
	functional (see~\cite{ACGT3} Section 2.3). 

	The key point is that $\mathcal{M}$ is decreasing under the
	flow (this is well-known for the modified Mabuchi functional,
	since the Calabi flow is its gradient flow). This
	follows from the computation
	\[ \begin{aligned} \frac{d\mathcal{M}(\phi_t)}{dt} &= \int_0^m (-\Psi
		S(\phi_t)'' + \phi_t S(\phi_t)'')(1+\tau)\,d\tau \\
		&= \int_0^m (\phi_t-\Psi)(S(\phi_t)-S(\Psi))''(1+\tau)\, d\tau +
		\int _0^m \phi_t S(\Psi)'' (1+\tau)\, d\tau \\
		&\leqslant -2\int_0^m (S(\phi_t)-S(\Psi))^2 (1+\tau)\,d\tau,
	\end{aligned} \]
	where we have used that $\Psi S(\Psi)''=0$ and $S(\Psi)''$ is a negative
	distribution. 
	
	On the other hand we have that
	\[ \mathcal{M}(\phi) \geqslant \int_0^m\log\phi\cdot
	(1+\tau)d\tau \geqslant -C_1\int_0^m \log\frac{\Theta}{\phi}\,
	d\tau -C_2,\]
	where $\Theta$ is a fixed momentum profile and $C_1,C_2$ are
	constants. Since $\log$ is
	concave we obtain
	\[ \mathcal{M}(\phi) \geqslant -C_3\log\int_0^m
	\frac{\Theta}{\phi}\,d\tau - C_4,\]
	for some constants $C_3,C_4$. The lemma that follows now implies
	that along the flow
	\[ \mathcal{M}(\phi_t)\geqslant -C\log(1+t) - D. \]
	Since $\mathcal{M}(\phi_t)$ is decreasing, we necessarily
	have that along a subsequence its derivative tends to zero, ie. 
	$S(\phi_t)\to S(\Psi)$ in $L^2$ (integrating with respect to
	$(1+\tau)d\tau$
	as usual). 
	Since $\Vert S(\phi_t)\Vert_{L^2}$ is
	decreasing along the flow, it follows that 
	\begin{equation}\label{eq:calabiminimise}
	  \lim_{t\to\infty}\Vert S(\phi_t)\Vert_{L^2}= \Vert S(\Psi)\Vert_{L^2}.
	\end{equation}
	
	Let us now take any
	subsequence $\phi_i$. Because of the uniform $H^2$-bound there is a
	subsequence also denoted by $\phi_i$ which converges weakly in $H^2$ to
	some limit. Now Equation~\ref{eq:calabiminimise} implies the convergence
	of the $H^2$-norms, which together with the weak convergence implies
	strong convergence in $H^2$. The limit then has to be $\Psi$ since the
	minimiser of the Calabi functional is unique
	(Proposition~\ref{prop:minimumprofile}). 
\end{proof}

\begin{lem}\label{lem:intbound}
  Let $\Theta : [0,m]\to\mathbf{R}$ be a
  momentum profile. 
  For the solution $\phi_t$ to the Calabi flow we have
  \begin{equation}\label{eq:intbound}
    \int_0^m\frac{\Theta}{\phi_t}\,d\tau < C(1+t)
  \end{equation}
  for some constant $C$. 
\end{lem}
\begin{proof}
  Let us define
  the functional
  \[ \mathcal{F}(\psi) =  \int_0^m \frac{\Theta}{\psi} - \log\frac{
  \Theta}{\psi}\,d\tau \]
  for any momentum profile $\psi$. Along the Calabi flow we have
  \[ \begin{aligned}
    \frac{d}{dt}\mathcal{F}(\phi_t) & = \int_0^m (\phi_t - \Theta)
    S(\phi_t)''\,d\tau 
     = \int_0^m (\phi_t - \Theta)''\, S(\phi_t)\,d\tau \\
     &\leqslant \left(\int_0^m (\phi_t''-\Theta'')^2\,d\tau\right) ^{1/2}
    (Cal(\phi_t)+C)^{1/2}. \end{aligned}
  \]
  The uniform $H^2$ bound on $\phi_t$ now implies that
  $\mathcal{F}(\phi_t) \leqslant C(1+t)$ 
  for some $C>0$. The result follows from the inequality $x-\log x >
  x/2$.
\end{proof}

\begin{rem}
	Note that when $m\leqslant k_2(k_2+2)$ the functional
	$\mathcal{M}$ is bounded below on the set of momentum profiles. This is
	because we can write
	\[ \mathcal{M}(\phi) = \int_0^m \left(\frac{\Psi}{\phi} -
	\log\frac{\Psi}{\phi}\right)\, (1+\tau)d\tau + \int_0^m \log\Psi\cdot
	(1+\tau)d\tau.\]
	Since $\Psi$ only vanishes at isolated points and to finite order, the
	integral of $\log\Psi$ is finite, so the inequality $\log x<x$ implies
	\[ \mathcal{M}(\phi) \geqslant \int_0^m \log\Psi\cdot (1+\tau) d\tau.\]
	In the case $m > k_2(k_2+2)$ however $\mathcal{M}$ is not bounded from
	below since now $\Psi$ vanishes on an interval. In particular as
	$\phi\to \Psi$, it is clear that $\mathcal{M}(\phi)\to -\infty$. 
\end{rem}

\section{Long time existence}\label{sec:longtimeexist}

The existence of the Calabi flow for a short time has been proved by
Chen-He~\cite{CH06} (also Guan~\cite{Gua05} for ruled manifolds). 
In the case when an extremal metric exists, the long time existence has
also been shown in~\cite{Gua05} for ruled manifolds. 

To show that the flow exists for all time 
we first need to show that $\phi_t(x)$ does not become zero in finite
time for $x\in(0,m)$. Let $\Theta$ be a fixed momentum profile, ie. a
non-negative function on $[0,m]$, strictly positive on the interior, and
satisfying the usual boundary conditions. We want to show

\begin{prop}\label{prop:c0bound}
	If $\phi_t$ is the solution to the Calabi flow, then
	$\sup\frac{\Theta(x)}{\phi_t(x)}$ does not blow up in finite time. 
\end{prop}
\begin{proof} This follows from Lemma~\ref{lem:intbound} and the
  following lemma. 
\end{proof}

\begin{lem} Given a constant $C>0$ there exists a constant $D>0$ such
  that if for a momentum profile $\psi$ we have 
  \[\int_0^m\frac{\Theta}{\psi}\,d\tau < C \quad \text{ and }\quad
  \Vert\psi\Vert_{C^{1,1/2}}<C,\]
  then
  \[\sup \Theta/\psi < D.\] 
\end{lem}
\begin{proof}
  Let us derive the estimate near the boundary first. 
  Because of the $C^{1,1/2}$ bound on $\psi$, there exists
  a constant $C_1$ such that
  \[ |\psi'(x)-\psi'(0)| < C_1\sqrt{x}, \]
  ie. 
  \[ \psi'(x) > 2-C_1\sqrt{x}. \]
  This implies that
  \[ \psi(x) > x\left( 2-\frac{2}{3}C_1\sqrt{x} \right), \]
  so that for $x < (3/2C_1)^2$ we have $\psi(x)>x$. 
  We can apply the same argument around $x=m$ as well, so we obtain a
  small constant $\delta$ such that
  \[ \mbox{if } x<\delta \mbox{ or } x > m-\delta, \mbox{ then }
  \frac{\Theta(x)}{\psi(x)} < D. \]

  Now we concentrate on the set $(\delta,m-\delta)$. On this set
  we have a uniform lower bound $\Theta(x)>\epsilon>0$ so we just need a
  lower bound on $\psi$. 
  There is a constant $C_2$
  such that $|\psi'(x)|<C_2$ for all $x$. Suppose that for 
  some $x\in(\delta,m-\delta)$ we have
  $\psi(x)<\epsilon/k$ where $k$ is large. Assume for simplicity that
  $x < m/2$. Then for $y<m/2-\delta$ we have
  \[ \psi(x+y) < \frac{\epsilon}{k} + C_2y. \]
  Writing $a=m/2-\delta$, this implies that 
  \[ 
    C > \int_0^m \frac{\Theta}{\psi}\,d\tau >
    \epsilon\int_0^a\frac{1}{\frac{\epsilon}{k} + C_2y}\,dy 
    > \frac{\epsilon}{C_2}\left[\log C_2a -
    \log\frac{\epsilon}{k}\right].
  \]
  Since this tends to infinity as $k\to\infty$, we get the required
  lower bound on $\psi(x)$ for $x\in(\delta,m-\delta)$. Combining this
  with the boundary estimate we obtain the statement of the
  lemma.
\end{proof}

Next we would like to estimate the derivatives of $\phi$ following the
calculation in Guan~\cite{Gua05}. Let us introduce the
functional 
\[ L(\phi) = \int_0^m (\phi S(\phi)'')^2\, (1+\tau)d\tau. \]
We want to show
\begin{lem}[Guan~\cite{Gua05}] 
	For $\phi_t$ a solution of the Calabi flow we have that
	$L(\phi_t)\leqslant C(t)$ for some function $C(t)$ defined for all $t$. 
\end{lem}
\begin{proof}
	All our constants will depend on $t$
	but will be finite for all $t$.
	All the integral norms will be with respect to the measure $d\tau$ and not
	$(1+\tau)d\tau$ as before. 

	In the proof we will repeatedly use the Hardy-type inequality
	\[ \Vert f \Vert_{L^2(0,m)} 
	\leqslant C \Vert \phi_t^{-k+1} (\phi_t^k f)'\Vert_{L^2(0,m)}
	\]
	for $k\geqslant 1$ and any $f\in C^1[0,m]$ with the constant $C$
	depending on $t$. 
	Using Proposition~\ref{prop:c0bound}, 
	this is easy to derive from the inequality
	\[ \int_{-1}^1 f(x)^2\, dx \leqslant C\int_{-1}^1 \left[
	(1-x^2)^{-k+1} ( (1-x^2)^k f(x))'
	\right]^2 \, dx. \]
	This in turn follows from the inequality
	\[ \int_0^1 f(x)^2\, dx \leqslant C\int_0^1 \left( x^{-k+1} (x^k
	f)'\right)^2\, dx \]
	for $f$ with $f(1)=0$,
	applied to the intervals $[-1,0]$ and $[0,1]$ separately
	(see~\cite{Gua05}). 

\begin{comment}
	We also use inequalities of the type 
	\[ \begin{aligned}
		\Vert \phi_t f'\Vert_{L^2} &\leqslant C\Vert (\phi_t f)'\Vert_{L^2},
		\\
		\Vert \phi_t^2 f''\Vert_{L^2} &\leqslant C\Vert (\phi_t^2 f)''
		\Vert_{L^2},
	\end{aligned} \]
	where again $f$ vanishes on the boundary. 
	These follow from the fact that we have a uniform $C^1$ bound on
	$\phi_t$ along the flow. For example the first inequality follows from
	\[ \Vert \phi_t f'\Vert_{L^2}\leqslant \Vert (\phi_t f)'\Vert_{L^2} + \Vert
	\phi_t' f\Vert_{L^2}\leqslant  \Vert (\phi_t f)'\Vert_{L^2} + C\Vert
	f\Vert_{L^2}, \]
	and
	\[ \Vert f\Vert_{L^2}= \Vert \phi_t^{-1}(\phi_t f)\Vert_{L^2}
	\leqslant C \Vert (\phi_t f)'\Vert_{L^2}. \]
\end{comment}

	Let us compute the derivative of $L(\phi_t)$. 
\[ \begin{aligned}
	\frac{d}{dt} L(\phi_t) &= 2\int_0^m (\phi_t S(\phi_t)'')^3\, (1+\tau)d\tau -
	\int_0^m \left[(1+\tau)\phi_t^2 S(\phi_t)''\right]''^2\frac{d\tau}{1+\tau}
	\\	
%	\int_0^m 2\phi_t S(\phi_t)'' \left
%	(\phi_t^2(S(\phi_t)'')^2 - \phi_t\left(\frac{1}{2(1+\tau)}[(1+\tau) \phi_t^2
%	S(\phi_t)'']\right)''\right)\, (1+\tau)d\tau \\
	&\leqslant C_1\int_0^m (\phi_t S(\phi_t)'')^3\,d\tau - C_2 \left\Vert
	\left(\phi_t^2 S(\phi_t)''\right)''\right\Vert^2_{L^2}.
\end{aligned} \]
Let us estimate the cubed term. We have
\[ \begin{aligned}
	\int_0^m (\phi_t S(\phi_t)'')^3\,d\tau &\leqslant  
	C_3\Vert \phi_t S(\phi_t)''\Vert_{C^0}
	L(\phi_t) \\
	&\leqslant C_4\Vert (\phi_t S(\phi_t)'')'\Vert_{L^2} L(\phi_t) \\
	&\leqslant C(\epsilon) L(\phi_t)^2 + \epsilon\Vert (\phi_t
	S(\phi_t)'')'\Vert^2_{L^2},
\end{aligned} \]
for any $\epsilon>0$ using Young's inequality. 
Using the uniform $H^2$-bound on $\phi_t$ and the Hardy-type inequality twice we
obtain
\[ \begin{aligned}
  \Vert (\phi_t^{-1}\cdot \phi_t^2
  S(\phi_t)'')'\Vert_{L^2} &\leqslant \Vert \phi_t^{-1}
  (\phi_t^2 S(\phi)'')'\Vert_{L^2} + \Vert \phi_t'
  S(\phi_t)''\Vert_{L^2}
  \\
  &\leqslant C_5 \Vert \phi_t^{-1} (\phi_t^2
  S(\phi_t)'')'\Vert_{L^2} \\
  &\leqslant C_6 \Vert (\phi_t^2 S(\phi_t)'')''\Vert_{L^2},
\end{aligned} \]
%  C\left\Vert \phi_t^2\frac{d^4
%S(\phi_t)}{d\tau^4}\right\Vert^2_{L^2}, \]
so if we choose $\epsilon$ small enough (depending on $t$), then we obtain the
inequality
\[ \frac{d}{dt} L(\phi_t) \leqslant C_1(t) L(\phi_t)^2. \]
This implies that
\[ \frac{d}{dt}\log L(\phi_t) \leqslant C_1(t) L(\phi_t),\]
ie. for any $T>0$ we have
\[ \log L(\phi_T) \leqslant \log L(\phi_0) + \sup_{t\in [0,T]}C_1(t) \int_0^T
L(\phi_t)\, dt.\]
Now Lemma~\ref{lem:Calabidec} gives a bound on the integral of $L(\phi_t)$ since the
Calabi functional is non-negative, so the proof is complete. 
\end{proof}

Now we need to use the inequality
\[ \Vert f\Vert_{L^2}^2 \leqslant C( \Vert \phi f'\Vert_{L^2}^2 + f(m/2)^2) \]
for all $f\in C^1(0,m)$ which
can be proved in the same way as the Hardy-type inequalities we
used before. This implies that
\[ \Vert S(\phi_t)\Vert^2_{C^0}\leqslant C_1\Vert S(\phi_t)'\Vert^2_{L^2} \leqslant
C_2\big[\Vert \phi_t S(\phi_t)''\Vert^2_{L^2} + (S(\phi_t)'(m/2))^2\big]. \]
The bound on $\Vert \phi_t S(\phi_t)''\Vert_{L^2}$ gives a bound on
$|S(\phi_t)'(x)-S(\phi_t)'(m/2)|$ for $x$ inside the interval $\left(\frac{
m}{3}, \frac{2m}{3}\right)$. The bound on $\Vert S(\phi_t)\Vert_{L^2}$ (the
Calabi functional decreases along the flow) then gives an apriori bound on
$S(\phi_t)'(m/2)$. 
Therefore as long as $L(\phi_t)$ remains bounded, we have a $C^2$ bound on 
$\phi_t$ (depending on
$t$). To obtain estimates for the higher derivatives of $\phi_t$ we could either
continue with similar integral estimates 
in the manner of~\cite{Gua05} or we can note that a $C^2$ bound on
the momentum profile implies a uniform bound on the Ricci curvature. According
to Chen-He~\cite{CH06} the Calabi flow exists for all time as long as the Ricci
curvature remains uniformly bounded. 

\bibliographystyle{hplain} 
\bibliography{mybib}

\begin{thebibliography}{10}

\bibitem{ACGT3}
V.~Apostolov, D.~M.~J. Calderbank, P.~Gauduchon, and C.~W.
  T\o{}nnesen-Friedman.
\newblock Hamiltonian 2-forms in {K}\"ahler geometry {III}, extremal metrics
  and stability, arXiv:math.DG/0511118.

\bibitem{Aub78}
T.~Aubin.
\newblock {\'E}quations du type {M}onge-{A}mp\`ere sur les variet\'es
  k\"ahl\'eriennes compactes.
\newblock {\em Bull. Sci. Math. (2)}, 102(1):63--95, 1978.

\bibitem{Cal82}
E.~Calabi.
\newblock Extremal {K}\"ahler metrics.
\newblock In S.~T. Yau, editor, {\em Seminar on Differential Geometry}.
  Princeton, 1982.

\bibitem{CH06}
X.~Chen and Weiyong He.
\newblock On the {C}alabi flow, arXiv:math.DG/0603523.

\bibitem{Chen06}
X.~X. Chen.
\newblock {Space of K\"ahler metrics III--On the lower bound of the Calabi
  energy and geodesic distance}, arXiv:math.DG/0606228.

\bibitem{Don97}
S.~K. Donaldson.
\newblock Remarks on gauge theory, complex geometry and four-manifold topology.
\newblock In Atiyah and Iagolnitzer, editors, {\em Fields Medallists'
  Lectures}, pages 384--403. World Scientific, 1997.

\bibitem{Don02}
S.~K. Donaldson.
\newblock Scalar curvature and stability of toric varieties.
\newblock {\em J. Differential Geom.}, 62:289--349, 2002.

\bibitem{Don04_1}
S.~K. Donaldson.
\newblock Conjectures in {K}\"ahler geometry.
\newblock In {\em Strings and geometry}, volume~3 of {\em Clay Math. Proc.},
  pages 71--78. Amer. Math. Soc., Providence, RI, 2004.

\bibitem{Don05}
S.~K. Donaldson.
\newblock Lower bounds on the {C}alabi functional.
\newblock {\em J. Differential Geom.}, 70(3):453--472, 2005.

\bibitem{Fuj92}
A.~Fujiki.
\newblock Moduli space of polarized algebraic manifolds and {K}\"ahler metrics
  [translation of {S}\^ugaku {\bf 42} (1990), no. 3, 231--243;].
\newblock {\em Sugaku Expositions}, 5(2):173--191, 1992.

\bibitem{Gua05}
D.~Guan.
\newblock Extremal-solitons and exponential ${C}^\infty$ convergence of the
  modified {C}alabi flow on certain $\mathbf{CP}^1$ bundles.
\newblock {\em preprint}, 2005.

\bibitem{HS02}
A.~Hwang and M.~A. Singer.
\newblock A momentum construction for circle-invariant {K}\"ahler metrics.
\newblock {\em Trans. Amer. Math. Soc.}, 354(6):2285--2325, 2002.

\bibitem{Hwang95}
A.~D. Hwang.
\newblock On the {C}alabi energy of extremal {K}\"ahler metrics.
\newblock {\em Internat. J. Math.}, 6(6):825--830, 1995.

\bibitem{Mab04_1}
T.~Mabuchi.
\newblock Stability of extremal {K}\"ahler manifolds.
\newblock {\em Osaka J. Math.}, 41, 2004.

\bibitem{GSz04}
G.~Sz\'ekelyhidi.
\newblock Extremal metrics and ${K}$-stability, arXiv:math.AG/0410401.

\bibitem{GSzThesis}
G.~Sz\'ekelyhidi.
\newblock {\em Extremal metrics and ${K}$-stability}.
\newblock PhD thesis, Imperial College, London, 2006, math.DG/0611002.

\bibitem{Tian90}
G.~Tian.
\newblock On {C}alabi's conjecture for complex surfaces with positive first
  {C}hern class.
\newblock {\em Invent. Math.}, 101(1):101--172, 1990.

\bibitem{Tian97}
G.~Tian.
\newblock K\"ahler-{E}instein metrics with positive scalar curvature.
\newblock {\em Invent. math}, 137:1--37, 1997.

\bibitem{TF97}
C.~T\o{}nnesen-Friedman.
\newblock {\em Extremal {K}\"ahler metrics on ruled surfaces}.
\newblock PhD thesis, Odense University, 1997.

\bibitem{Yau78}
S.-T. Yau.
\newblock On the {R}icci curvature of a compact {K}\"ahler manifold and the
  complex {M}onge-{A}mp\`ere equation {I}.
\newblock {\em Comment. Pure Appl. Math.}, 31:339--411, 1978.

\end{thebibliography}

\end{document}